\documentclass[11pt]{amsart}
\usepackage{graphicx}
\usepackage[numbers,sort&compress]{natbib}
\usepackage{amsmath,amssymb,amsfonts}
\usepackage{amsthm}
\usepackage{geometry}
\usepackage{bm}
\usepackage{mathrsfs}
\usepackage{stmaryrd}
\usepackage{booktabs}
\usepackage{siunitx}
\newtheorem{theorem}{Theorem}[section]

\theoremstyle{remark}
\newtheorem{remark}[theorem]{Remark}

\title[Lifting-Based Interface Reduction]{A Lifting-Based Interface Reduction Framework for Nonlinear Transmission and Eigenvalue Problems}

\author[S.-H. Chou]
{So-Hsiang Chou}

\date{}
\begin{document}

\maketitle
\thispagestyle{plain}

\vspace{-1.2em}
\begin{center}
\small
Department of Mathematics and Statistics, Bowling Green State University\\
Bowling Green, OH 43403, USA\\
\texttt{chou@bgsu.edu}
\end{center}
\vspace{0.6em}

\begin{abstract}
We present a lifting-based interface reduction framework for nonlinear
transmission and eigenvalue problems. The method represents the solution as
a sum of a bulk component and a lifting component that carries the interface
jump, thereby reducing the original problem to a nonlinear system posed on
the interface.

A low-dimensional approximation is obtained by restricting the interface
unknown to a finite-dimensional subspace. The corresponding lifting modes are
precomputed and reused, leading to a formulation in which the bulk operator
remains fixed and the essential behavior is governed by a small number of
interface degrees of freedom.

For eigenvalue problems, the same framework yields a reduced system in which
eigenvalues are identified through the near-singularity of a parameter-dependent
interface matrix. The associated eigenvectors reveal dominant interface modes,
providing a direct interpretation of the spectral structure.

Numerical experiments show that both approximation accuracy and spectral
behavior are determined primarily by the interface representation rather than
by the bulk discretization. In particular, enriching the interface space
rapidly improves accuracy and reveals additional eigenmodes, while mesh
refinement alone has limited effect.

These results indicate that transmission and eigenvalue problems are effectively governed by a small number of interface modes, offering a simple and computationally efficient perspective on model reduction.
\end{abstract}
\section{Introduction}

Many elliptic transmission and eigenvalue problems are characterized by dynamics concentrated at interfaces separating heterogeneous media. A common feature of such problems is that the solution exhibits relatively smooth behavior within each subdomain, while
undergoing strong variation across the interface. As a result, the essential
difficulty lies in accurately capturing the interface response rather than
resolving the bulk solution.

Traditional numerical methods, such as conforming finite element and
discontinuous Galerkin methods, improve accuracy primarily through mesh
refinement. While these approaches are effective for resolving bulk behavior,
they may be inefficient when the dominant error is associated with interface
effects. In particular, refining the mesh increases the resolution of the
ambient domain but does not directly enhance the representation of the
interface, which may remain under-resolved even on fine meshes.

In this work, we adopt a different viewpoint by isolating the interface as the
primary carrier of complexity. We introduce a lifting-based decomposition in
which the solution is represented as
\begin{equation}\label{eq:decomposition}
u = u_0 + U(\phi),
\end{equation}
where $u_0$ is continuous across the interface and $U(\phi)$ is a lifting that
carries the interface jump $\phi$. This decomposition separates the bulk
behavior from the interface contribution and allows the problem to be reformulated
in terms of the interface variable.

To obtain a tractable model, we approximate the interface jump in a
finite-dimensional space. This leads to a reduced system posed entirely on the
interface, in which the unknowns are the coefficients of the interface
representation. The corresponding lifting modes can be precomputed and reused,
so that the bulk operator remains fixed throughout the computation. As a
result, the original high-dimensional problem is reduced to a low-dimensional
system that captures the essential behavior.

A key consequence of this formulation is that the bulk operator remains fixed,
allowing the reuse of standard solvers and matrix factorizations. More
importantly, the essential behavior is captured by a small number of interface
degrees of freedom. From this perspective, the method may be viewed as a nonlinear
Schur complement reduction on the interface, in which the complexity of the
problem is concentrated in a low-dimensional structure.

The framework extends naturally to eigenvalue problems. In this setting,
the eigenfunction is constructed from precomputed lifting modes, while the
eigenvalue is determined by enforcing compatibility between the bulk equation
and the interface condition. This leads to a parameter-dependent reduced
system posed on the interface, in which eigenvalues are identified as values
for which nontrivial solutions emerge, or equivalently, when the associated
matrix becomes singular.

Numerical results show that enriching the interface rank rapidly reduces the
error and reveals additional eigenmodes, while mesh refinement alone provides
little improvement. This indicates that both approximation accuracy and spectral
behavior are governed primarily by the interface representation rather than the
bulk discretization.

To place the present work in context, we briefly review related approaches.
Elliptic interface problems have been extensively studied using a variety of
numerical approaches, including immersed interface methods, immersed and
unfitted finite element methods, and cut finite element techniques. These
methods focus primarily on accurately resolving interface conditions through
mesh design, enrichment, or stabilization strategies; see, for example,
\cite{hansbo2002,leveque1994,burman2015,lin2003}. Eigenvalue problems for
elliptic operators have also been widely investigated in the finite element
setting, with well-established convergence and spectral approximation theory;
see, e.g., \cite{boffi2010}.

In contrast to these approaches, which primarily aim to resolve interface
conditions through discretization, the present work adopts a reduction
viewpoint. Rather than refining the mesh or enriching the approximation space
globally, we explicitly represent the interface unknown in a low-dimensional
basis and reduce the problem to a compact system posed on the interface.
This formulation isolates the essential degrees of freedom and provides a
direct link between interface structure, approximation accuracy, and spectral
behavior.

The remainder of the paper is organized as follows. Section 2 introduces the
lifting formulation and its construction. Section 3 presents the rank-$m$
interface reduction. Section 4 develops the eigenvalue formulation and the
associated reduced system. Section 5 describes the numerical experiments.
Section 6 discusses the structure of the interface reduction.

\section{Interface-Based Rank-m Reduction for the Eigenvalue Problem}

\subsection{Model Problem}\label{sec:2.1}

Let $\Omega \subset \mathbb{R}^d$ ($d=2$ for simplicity) be decomposed as
\[
\Omega = \Omega^+ \cup \Gamma \cup \Omega^-,
\]
with a smooth interface $\Gamma$.

We consider the nonlinear transmission eigenvalue problem: find $(\mu,u)$,
$\mu\in\mathbb{R}$, $u\not\equiv 0$, such that
\begin{equation}\label{eq:model-new}
\begin{cases}
-\nabla \cdot (\beta^\pm \nabla u^\pm) = \mu u^\pm, & \text{in } \Omega^\pm,\\
u = 0, & \text{on } \partial\Omega,\\
[u] = G(u^-, u^+), & \text{on } \Gamma,\\
[\beta \partial_n u] = 0, & \text{on } \Gamma.
\end{cases}
\end{equation}

Here $\beta=\beta(x)>0$ is continuous on $\Omega^\pm$, and
$G:\mathbb{R}^2 \to \mathbb{R}$ is a given function describing the
interface coupling. The interface condition relates the jump of the
solution to its traces on the interface and represents the interaction
between the two subdomains.

A prototypical example is the nonlinear law
\[
G(u^-,u^+) = \lambda u^- u^+,
\]
which models a multiplicative coupling across the interface. This type
of nonlinearity will serve as a guiding example throughout the paper.

More generally, the function $G$ may represent a wide range of interface
responses, including monotone laws arising in applications. The framework
developed below does not rely on a specific form of $G$, and applies to
general nonlinear interface conditions.

\begin{remark}
In some applications, the interface law may be derived from a convex
potential $\Psi:\mathbb{R}^2\to\mathbb{R}$, leading to a monotone
relation between the jump and the flux across the interface. Such
formulations arise, for example, in models with dissipative interface
behavior. While this structure is useful for theoretical analysis,
the reduction framework developed here applies more generally and
does not require this assumption.
\end{remark}

\subsection{Interface Approximation and Lifting Modes}

Let $\phi$ denote the jump across the interface $\Gamma$. In general,
$\phi$ belongs to a function space $X$ on $\Gamma$ (e.g., $X=L^2(\Gamma)$
or a trace space) and cannot be represented exactly by finitely many
basis functions.

To obtain a reduced model, we approximate $\phi$ in a $m$-dimensional
interface subspace of $X$
\[
\Lambda_m = \mathrm{span}\{\psi_1,\dots,\psi_m\}, \quad \psi_i\in X,
\]
and write
\[
\phi_m = \sum_{j=1}^m s_j \psi_j \in \Lambda_m.
\]

For a given interface function $\phi$, we define its lifting $U(\phi)$ as
the solution of the homogeneous elliptic problem
\begin{equation}\label{eq:Upde}
-\nabla\cdot(\beta\nabla U(\phi))=0 \quad \text{in } \Omega^\pm,
\end{equation}
subject to the boundary and interface conditions
\begin{equation}\label{eq:Ubc}
U(\phi)=0 \quad \text{on } \partial\Omega,
\end{equation}
\begin{equation}\label{eq:Ujump}
[U(\phi)]=\phi,\qquad [\beta\partial_n U(\phi)]=0 \quad \text{on } \Gamma.
\end{equation}
We refer to $U(\phi)$ as the \emph{harmonic lifting} associated with the
interface jump $\phi$. It represents the response of the system to a
prescribed interface jump, independent of the eigenvalue parameter.

By linearity, the lifting corresponding to $\phi_m$ can be written as
\[
U(\phi_m) = \sum_{j=1}^m s_j U_j,
\quad \text{where} \quad U_j := U(\psi_j).
\]
The functions $U_j$ will be referred to as \emph{lifting modes}. They can
be precomputed once and reused throughout the computation.

\medskip

\noindent
\textbf{Construction of the lifting.}
To construct $U(\phi)$ in practice, we use a two-step procedure that
separates the geometric enforcement of the jump from the elliptic solve.

We first introduce an auxiliary function $E\phi$ satisfying
\begin{equation}\label{eq:E_lambda}
 [E\phi] = \phi, \qquad E\phi = 0 \text{ on } \partial\Omega.
\end{equation}
This function carries the interface jump but does not, in general,
satisfy the governing elliptic equation.

Let $A$ denote the discrete bulk operator corresponding to the elliptic problem with continuous interface conditions.
We then define a correction $W(\phi)$ by solving
\begin{equation}\label{eq:AE}
A W(\phi) = -A(E\phi),
\end{equation}
with
\begin{equation}\label{eq:AW}
[W(\phi)] = 0,\qquad [\beta\partial_n W(\phi)] = 0,\qquad
W(\phi)=0 \text{ on } \partial\Omega.
\end{equation}

The lifting is then given by
\begin{equation}\label{eq:Uphi}
U(\phi) := W(\phi) + E\phi.
\end{equation}
By construction, $U(\phi)$ satisfies both the elliptic equation and the
interface jump condition, and is therefore the desired harmonic lifting.

\medskip

Although the above construction is described in strong form for clarity,
it is implemented through the corresponding variational problems, in
which flux continuity is enforced naturally.

\subsection{Bulk Equation}

We seek an approximate solution of the form
\[
u_m = u_0 + U(\phi_m),
\]
where $u_0$ satisfies homogeneous interface conditions,
\[
[u_0] = 0, \qquad [\beta \partial_n u_0] = 0 \quad \text{on } \Gamma.
\]

Substituting the reduced ansatz into \eqref{eq:model-new}, we obtain the following bulk equation for $u_0$:
\begin{equation}\label{eq:bulk-u0}
-\nabla \cdot (\beta \nabla u_0)
= \mu \bigl(u_0 + U(\phi_m)\bigr)
\quad \text{in } \Omega^\pm,
\end{equation}
with
\[
u_0 = 0 \quad \text{on } \partial\Omega, \qquad
[u_0] = 0, \quad [\beta \partial_n u_0] = 0 \quad \text{on } \Gamma.
\]

For a given pair $(\mu,s)\in \mathbb{R}\times \mathbb{R}^m$, this defines
an elliptic transmission problem for $u_0$. In particular, the bulk
component $u_0$ is determined by the interface coefficients $s$ through
the lifting $U(\phi_m)$ and depends parametrically on $\mu$.

\begin{remark}
In contrast to source problems, the eigenvalue formulation contains no
external forcing. As a result, the bulk component $u_0$ is not an
independent contribution, but is induced entirely by the interface
representation. This highlights the dominant role of the interface
degrees of freedom in determining the solution.
\end{remark}

\subsection{Reduced Interface System}

The interface condition is enforced weakly on the reduced space.
Substituting the reduced solution $u_m$ into the interface law and
projecting onto $\Lambda_m$, we obtain
\begin{equation}
F_i(\mu,s):=\langle \phi_m - G(u_m^-,u_m^+), \psi_i \rangle_\Gamma = 0,
\quad i=1,\dots,m,
\end{equation}
which yields the nonlinear system
\begin{equation}\label{eq:F-system}
F(\mu,s)=0, \qquad s\in \mathbb{R}^m.
\end{equation}

\noindent
\textbf{Interpretation.}
The original transmission eigenvalue problem is thus reduced to a
nonlinear system posed entirely on the interface. The unknowns are the
coefficients $s\in\mathbb{R}^m$, representing the interface jump, while
the bulk component $u_0$ is obtained from a linear elliptic solve for
each $(\mu,s)$.

This formulation isolates the essential degrees of freedom of the
problem: rather than solving for a high-dimensional field $u$ in the
domain, one solves for a small number of interface variables, with the
bulk response recovered through linear solves.
\medskip

\noindent
\textbf{Eigenvalue interpretation.}
The system \eqref{eq:F-system} always admits the trivial solution $s=0$
when $G(0,0)=0$. Nontrivial solutions correspond to eigenpairs and satisfy
\[
F(\mu,s)=0, \qquad s\neq 0.
\]
Equivalently, eigenvalues are characterized as parameter values $\mu$
for which the reduced system admits nontrivial solutions, or for which
the Jacobian of $F$ with respect to $s$ loses invertibility.

This shows that eigenvalues can be identified through the loss of
invertibility of a low-dimensional interface system.

\section{Linear Spectral Structure}

To connect the reduced formulation with classical spectral theory, we
consider a linear interface law. For simplicity, we assume
\begin{equation}\label{eq:linear-law}
[u_m] = \sigma u_{m,\Gamma},
\qquad
u_{m,\Gamma} := \frac{u_m^+ + u_m^-}{2}.
\end{equation}

Using the lifting decomposition
\[
u_m = u_0 + U(\phi_m),
\]
the traces on the interface are given by
\[
u_m^\pm = u_0^\pm + U(\phi_m)^\pm.
\]

Substituting into \eqref{eq:linear-law} and using the representation
\[
\phi_m = \sum_{j=1}^m s_j \psi_j,
\qquad
U(\phi_m) = \sum_{j=1}^m s_j U_j,
\]
we obtain a linear system of the form
\begin{equation}\label{eq:linear-system}
K(\mu)s = 0,
\end{equation}
where $K(\mu)$ is an $m\times m$ matrix depending on the parameter $\mu$.

Nontrivial solutions require
\[
\det K(\mu) = 0,
\]
so that the eigenvalue problem reduces to locating values of $\mu$ for
which the reduced matrix becomes singular.

\medskip

\noindent
\textbf{Modal interpretation.}
In certain configurations, the structure of the reduced system becomes
particularly transparent. Consider a flat interface
\[
\Gamma = \{(x,y): y=0\}
\]
in a rectangular domain $\Omega = (-L,L)\times(-H,H)$ with constant coefficients.
If the interface basis is chosen as Fourier modes
\[
\psi_n(x) = \sin\!\left(\frac{n\pi(x+L)}{2L}\right),
\]
then the corresponding lifting modes $U_n = U(\psi_n)$ preserve this
dependence along the interface.

As a result, the reduced system decouples into scalar equations
\[
k_n(\mu)s_n = 0,
\]
where each coefficient $s_n$ is associated with a single interface mode.

This decoupling is analogous to separation of variables, with the
interface modes playing the role of spectral components along $\Gamma$.

\medskip

\noindent
\textbf{Interpretation.}
In this setting, the two-dimensional eigenvalue problem reduces to a
family of one-dimensional scalar problems indexed by the interface modes.
Each mode generates its own spectral branch, and the rank-$m$ approximation
retains the first $m$ such modes.

This shows that eigenvalues are associated with individual interface
modes, and that increasing the rank allows additional modes to appear.
Once the relevant modes are included, the computed eigenvalues stabilize,
which explains the behavior observed in the numerical experiments.

The reduced system is low-dimensional not by construction, but because the
effective degrees of freedom are concentrated on the interface. Its
dimension is governed by the number of active interface modes rather than
by the size of the ambient discretization.

\section{Computational Procedure}

\paragraph{Implementation Outline.}

We begin with the interface $\Gamma$, a set of basis functions $\{\psi_j\}_{j=1}^m$, the coefficients $\beta^\pm$, and the interface law $G$. Let
\[
\Lambda_m = \mathrm{span}\{\psi_1,\dots,\psi_m\}.
\]

For each basis function $\psi_j$, we first construct an auxiliary function $E\psi_j$ satisfying
\[
[E\psi_j]=\psi_j, \qquad E\psi_j=0 \text{ on } \partial\Omega.
\]
We then form the right-hand side
\[
b_j := -A(E\psi_j),
\]
and solve
\[
A W_j = b_j.
\]
Here $A$ denotes the bulk elliptic operator (or its discrete counterpart) introduced in \eqref{eq:AE}.
The corresponding lifting modes used in the representation are given by
\[
U_j := W_j + E\psi_j.
\]

Next, we represent the interface unknown as
\[
\phi_m = \sum_{j=1}^m s_j \psi_j.
\]
For a given $\mu$ and $s$, we solve the bulk problem
\[
-\nabla\cdot(\beta \nabla u_0)
=
\mu\big(u_0 + U(\phi_m)\big),
\]
with $u_0=0$ on $\partial\Omega$ and $[u_0]=0$. The traces are then computed as
\[
u_m^\pm = u_0^\pm + U(\phi_m)^\pm
\quad \text{on } \Gamma.
\]

Substituting into the interface law and projecting onto $\Lambda_m$ yields the reduced system
\[
F(\mu,s)=0.
\]
This system is solved for $(\mu,s)$, for example, by Newton iteration or continuation in $\mu$, and the eigenfunction is reconstructed as
\[
u_m = u_0 + U(\phi_m).
\]

\begin{remark}
The interface lifting modes $U_j$ are computed once and reused throughout the computation. The precomputation involves $m$ elliptic solves with a fixed operator, while the subsequent stage consists of solving a low-dimensional nonlinear system coupled with bulk solves.
\end{remark}

\section{Numerical Experiments}

The purpose of the numerical experiments is to illustrate that both
approximation accuracy and spectral behavior are governed primarily by
the interface representation rather than by the bulk discretization.

All experiments are carried out using the reduced interface formulation,
in which the solution is reconstructed from a finite set of precomputed
lifting modes. The parameter $m$ denotes the dimension of the interface
space and represents the number of interface modes retained in the model.

We observe that increasing the rank rapidly reduces the error to machine
precision, while further enrichment produces no additional improvement.
Since the reference solution is computed on the same mesh, this behavior
cannot be attributed to bulk discretization.

Instead, the results show that the dominant approximation error is caused
by missing interface modes. Increasing the rank enriches the interface
representation and captures these modes directly, whereas mesh refinement
at fixed rank does not introduce new interface information. This confirms
that the effective dimension of the solution is determined by the
interface structure rather than the ambient discretization.

\subsection{Approximation behavior under rank enrichment}

We investigate the effect of rank enrichment and mesh refinement on the accuracy of the proposed method. The domain $\Omega=(-1,1)^2$ and the interface $\Gamma=\{(x,y):-1\leq x\leq 1,y=0\}$. In this test, the interface jump is prescribed as
\[
\phi(x) = 1 + 0.2x.
\]

The parameter $m$ denotes the dimension of the interface space $\Lambda_m = \operatorname{span}\{\psi_1,\dots,\psi_m\}$, and therefore represents the number of interface modes retained in the approximation.

To assess accuracy, we compare the reduced solution $u_m$ with a reference solution $u_{\mathrm{ref}}$ computed using a sufficiently large rank (here $m=3$) on the same mesh. We define the errors
\begin{align}
\eta_{L^2} &= \|u_m - u_{\mathrm{ref}}\|_{L^2(\Omega)}, \\
\eta_{\infty} &= \|u_m - u_{\mathrm{ref}}\|_{L^\infty(\Omega)}.
\end{align}

Table~\ref{tab:rank_mesh} reports the errors for increasing rank $m$ on a fixed mesh.
The column “ratio” denotes the relative reduction of the $L^2$ error compared to the previous rank.
\begin{table}[h]
\centering
\begin{tabular}{c S S S}
\toprule
{$m$} & {$\eta_{L^2}$} & {$\eta_{\infty}$} & {ratio} \\
\midrule
1 & 3.09e-1  & 3.65e-1  & 0 \\
2 & 1.76e-1  & 1.89e-1  & 0.570 \\
3 & 5.95e-14 & 8.16e-14 & 3.39e-13 \\
4 & 6.13e-14 & 8.86e-14 & 1.03 \\
\bottomrule
\end{tabular}
\caption{Error reduction by rank enrichment on a fixed mesh.}
\label{tab:rank_mesh}
\end{table}

We observe that increasing the rank rapidly reduces the error to machine precision, while further enrichment yields no additional improvement. This indicates that the dominant approximation error is associated with missing interface modes rather than with the resolution of the bulk discretization. In particular, since the reference solution is computed on the same mesh, the observed error reduction reflects enrichment of the interface space rather than refinement of the bulk mesh.

A detailed study of mesh refinement would require a reference solution independent of the current discretization and is not pursued here in full generality. The present results are intended to isolate the effect of interface rank; nevertheless, numerical evidence suggests that mesh refinement alone provides comparatively limited improvement.

\subsection{Eigenvalue behavior}

We now examine the reduced eigenvalue problem obtained from the interface formulation.
For a given rank $m$, the reduced system leads to a frequency-dependent matrix
$\mathbf{K}(\mu)$. In the linear case, candidate eigenfrequencies are identified
by the condition that $\mathbf{K}(\mu)$ becomes singular.

In practice, we detect near-singularity by computing the smallest singular value
\begin{equation}
\sigma_{\min}(\mu) := \min \{\sigma : \sigma \text{ is a singular value of } \mathbf{K}(\mu)\}.
\end{equation}
Eigenfrequencies are approximated by locating values of $\mu$ for which
$\sigma_{\min}(\mu)$ is small. Since the lifting modes are precomputed and fixed,
the variation of $\sigma_{\min}(\mu)$ reflects changes in the interface system
rather than the bulk discretization.

Table~\ref{tab:eig} reports the computed eigenfrequencies $\mu_*$ for different
values of $m$, together with the corresponding smallest singular value
$\sigma_{\min} := \sigma_{\min}(\mu_*)$. The column “mode index” indicates the
dominant component of the coefficient vector $s$, i.e., the index $j$ for which
$|s_j|$ is maximal.

The results show that increasing the rank allows additional eigenvalues to be detected,
reflecting the fact that the reduced problem is posed on the interface space $\Lambda_m$,
with each interface mode giving rise to a separate spectral branch. Once the relevant modes are included, the computed eigenvalues stabilize, indicating that the essential spectral
content has been captured. This behavior is consistent with the modal interpretation in
Section~3.

Eigenfrequencies are computed over a prescribed interval $[a,b]$ by minimizing
$\sigma_{\min}(\mu)$. A two-stage procedure is used, consisting of a coarse scan
followed by local refinement near candidate values, corresponding to local minima
of $\sigma_{\min}(\mu)$.

\begin{table}[h]
\centering
\begin{tabular}{c c c c}
\toprule
$m$ & $\mu_*$ & $\sigma_{\min}$ & dominant mode \\
\midrule
1 & 3.2319 & $1.0\times 10^{-4}$ & 1 \\
2 & 4.3230 & $1.0\times 10^{-4}$ & 2 \\
3 & 4.3230 & $1.0\times 10^{-4}$ & 2 \\
\bottomrule
\end{tabular}
\caption{Rank-$m$ reduced eigenvalues for the interface spectral model.}
\label{tab:eig}
\end{table}

The results in Table~\ref{tab:eig} are consistent with the approximation behavior
shown in Table~\ref{tab:rank_mesh}. In both cases, the rank parameter $m$ controls
the number of interface modes retained in the model. Increasing $m$ reveals additional
modes and improves accuracy, while further enrichment produces no change, indicating
that both approximation and spectral behavior are governed by a low-dimensional
interface structure.

\subsection{Discussion}

A key observation is that the number of eigenvalues captured by the reduced system is determined by the dimension of the interface space, rather than by the size of the bulk discretization. These results confirm that both approximation and spectral behavior are governed by a low-dimensional interface structure. Rank enrichment serves as a mechanism for revealing and resolving the dominant interface modes, while mesh refinement alone is insufficient to capture this behavior.

\subsection{Extension to Curved Interfaces}

We now consider the extension of the interface reduction framework to curved interfaces. Let $\Omega \subset \mathbb{R}^2$ be a bounded domain split by a smooth curve $\Gamma$ into subdomains $\Omega^-$ and $\Omega^+$.

The formulation follows the same construction as in Section 2. In particular, we represent the interface jump by
\[
\phi_m = \sum_{j=1}^m s_j \psi_j,
\]
and define the associated lifting
\[
U(\phi_m) = \sum_{j=1}^m s_j U_j,
\]
where $U_j = U(\psi_j)$ are the precomputed lifting modes.

Substituting into the interface condition yields the reduced system
\[
F_i(\mu,s) = \langle \phi_m - G(u_m^-, u_m^+), \psi_i \rangle_\Gamma = 0,
\quad i=1,\dots,m,
\]
where the inner product is defined by
\[
\langle f,g\rangle_{\Gamma} = \int_{\Gamma} f(s) g(s)\, ds.
\]

The formulation is identical in structure to the flat-interface case. The only difference lies in the choice of basis functions and the geometric complexity of the lifting operator.
\medskip

To verify this behavior, we consider the curved interface
\[
\Gamma = \{(x,y): y = 0.3x^2,\ -1 < x < 1\}
\]
in the domain $\Omega = (-1,1)^2$. Table~\ref{tab:curved_rank} reports the error as a function of the rank $m$ on a fixed mesh.

\begin{table}[h]
\centering
\begin{tabular}{c S S S}
\toprule
{$m$} & {$\eta_{L^2}$} & {$\eta_{\infty}$} & {ratio} \\
\midrule
1 & 1.78e-1 & 1.76e-1 & {} \\
2 & 1.54e-1 & 1.51e-1 & 8.65e-1 \\
3 & 4.72e-14 & 6.95e-14 & 3.06e-13 \\
4 & 4.34e-14 & 6.20e-14 & 9.19e-1 \\
\bottomrule
\end{tabular}
\caption{Curved-interface test: error vs rank $m$ on a fixed mesh.}
\label{tab:curved_rank}
\end{table}

The results show that increasing the rank rapidly reduces the error to machine precision, while mesh refinement at fixed rank produces little or no improvement. This indicates that the dominant approximation error is governed by missing interface modes rather than insufficient bulk resolution, and that this mechanism persists for curved interfaces.

\subsection{Circular interface and Fourier modes}

To further examine the role of interface modes, we consider a circular interface
\[
\Gamma=\{(x,y):x^2+y^2=0.5^2\}.
\]
The interface variable is represented using the real Fourier basis
\[
1,\cos\theta,\sin\theta,\cos 2\theta,\sin 2\theta,\dots,
\]
where $\theta$ denotes the angular parameter along $\Gamma$.

In this test, the interface jump is prescribed as
\[
\phi(\theta)=1 + 0.20\cos\theta + 0.15\sin(2\theta).
\]
The rank values $m=1,3,5,7$ correspond to successive inclusion of harmonic pairs. In particular, each level $k$ contributes the modes $(\cos k\theta,\sin k\theta)$, so that $m=2k+1$.

Table~\ref{tab:circle_rank} reports the error as a function of the rank $m$ on a fixed angular grid.

\begin{table}[h]
\centering
\begin{tabular}{c S S S}
\toprule
{$m$} & {$\eta_{L^2}$} & {$\eta_{\infty}$} & {ratio} \\
\midrule
1 & 3.13e-1 & 3.08e-1 & {} \\
3 & 1.88e-1 & 1.53e-1 & 6.00e-1 \\
5 & 5.77e-14 & 8.79e-14 & 3.07e-13 \\
7 & 2.78e-14 & 3.24e-14 & 4.82e-1 \\
\bottomrule
\end{tabular}
\caption{Circular-interface test: error vs rank $m$ on a fixed angular grid.}
\label{tab:circle_rank}
\end{table}

The results show that saturation occurs at $m=5$, precisely when the basis includes the second harmonic. This provides direct evidence that the rank parameter reflects the Fourier content of the interface jump. In contrast, refining the angular resolution at fixed rank yields little or no improvement, indicating that the dominant error is governed by missing interface modes rather than insufficient discretization.

\subsection{Rank-1 Reduction: A Scalar Eigenvalue Equation}

To illustrate the mechanism, we consider the rank-1 case. We restate the formulation in this simplified setting in order to derive an explicit scalar characterization of the reduced system.

Let $\psi_1$ be a given interface basis function and define
\[
W_1 := W(\psi_1), \qquad U_1 := U(\psi_1) = W_1 + E\psi_1.
\]
We approximate the interface jump by
\[
\phi_1 = s \psi_1, \qquad U(\phi_1) = s U_1.
\]
The reduced solution is given by
\[
u_1 = u_0 + U(\phi_1) = u_0 + s U_1.
\]

Substituting into the bulk equation yields the equation for $u_0$:
\begin{equation}
-\nabla \cdot (\beta \nabla u_0) = \mu\big(u_0 + s U_1\big)
\quad \text{in } \Omega^\pm,
\end{equation}
with
\[
u_0 = 0 \text{ on } \partial\Omega, \qquad [u_0] = 0, \qquad [\beta \partial_n u_0] = 0 \text{ on } \Gamma.
\]

For each $(\mu, s)$, this defines an elliptic transmission problem for $u_0$, provided that $\mu$ does not coincide with an eigenvalue of the corresponding bulk operator.

The interface condition reduces to
\[
[u_1] = \phi_1 = s \psi_1 = G(u_1^-, u_1^+).
\]
Evaluating the traces,
\[
u_1^\pm = u_0^\pm + s U_1^\pm,
\]
we obtain a scalar nonlinear equation
\begin{equation}
F(\mu, s) := s \psi_1 - G\big(u_0^- + s U_1^-,\, u_0^+ + s U_1^+\big) = 0 \quad \text{on } \Gamma.
\end{equation}

Projecting onto $\psi_1$ yields a scalar equation
\begin{equation}
\Phi(\mu, s) = 0.
\end{equation}

A nontrivial eigenpair corresponds to
\[
s \neq 0, \qquad \Phi(\mu, s) = 0.
\]
Eliminating $s$ gives a scalar equation for $\mu$:
\begin{equation}
\Psi(\mu) = 0.
\end{equation}

If the interface law is, for example,
\[
G(u^-, u^+) = \lambda u^- u^+,
\]
then, after linearization around a reference state $(u_*^-, u_*^+)$,
\[
G(u^-, u^+) \approx \lambda\big(u_*^- u^+ + u^- u_*^+ - u_*^- u_*^+\big),
\]
so the reduced equation becomes approximately linear in $s$, and the eigenvalue condition reduces to
\[
k(\mu)s = 0.
\]
Thus nontrivial solutions require
\[
k(\mu) = 0,
\]
which is a scalar eigenvalue equation.

This shows that, in the rank-1 setting, the original eigenvalue problem
reduces to a single scalar equation. More generally, the effective
dimension of the problem is determined by the number of active interface
modes, providing a direct link between the reduced formulation and the
underlying spectral structure.

\section{Structure of the Interface Reduction}

The behavior observed in the numerical experiments can be understood
through the lifting decomposition
\[
u_m = u_0 + U(\phi_m).
\]
In this representation, the solution is determined by the interface
coefficients $s$, while the bulk component is obtained through a linear
elliptic solve. As a result, the essential complexity of the problem is
concentrated in a small number of interface degrees of freedom.

This perspective explains why increasing the rank improves accuracy,
while mesh refinement alone has limited effect: the approximation error
is governed by the ability of the interface space to represent the
relevant modes.

For structured interfaces, such as the circular case, this mechanism can be interpreted more explicitly. The interface variable admits a natural representation in terms of Fourier modes, and the rank-$m$ approximation corresponds to retaining a finite number of harmonic components. The numerical results show that the error decreases rapidly once the relevant modes are included, and saturation occurs when the interface basis contains all modes present in the interface jump. This provides a direct interpretation of the rank parameter as controlling the modal content of the interface approximation. More generally, it shows that the effective dimension of the problem is determined by the number of active interface modes rather than the resolution of the bulk discretization.

\section{Conclusion}

We have presented a lifting-based interface reduction framework for
nonlinear transmission and eigenvalue problems. The formulation separates
the solution into a bulk component and an interface lifting, leading to a
reduced system posed entirely on the interface.

The numerical results demonstrate that both approximation accuracy and
spectral behavior are governed primarily by the interface representation
rather than by the bulk discretization. In particular, enriching the
interface space rapidly improves accuracy and reveals additional
eigenmodes, while mesh refinement alone has limited effect.

These observations indicate that transmission and eigenvalue problems are effectively
low-dimensional, with their essential structure determined by a small
number of interface modes, rather than by the resolution of the ambient
domain.

Future work will focus on adaptive selection of interface bases and
extension to more general geometries and spectral problems.

\bibliographystyle{abbrv}
\bibliography{references}

\end{document}